\documentclass[12pt,leqno]{article}
\usepackage{amstext,amsmath,amsthm,amsfonts,amssymb,amscd}
\usepackage{graphics}
\usepackage[latin1]{inputenc}

\usepackage{graphicx}
\usepackage{color}

\renewcommand{\thefootnote}

\begin{document}

\centerline{\LARGE\bf Rigidity of Convex Surfaces in }

\medskip

\centerline{\LARGE\bf Homogeneous Spaces}

\vskip .4in

\centerline{\bf by}

\vskip .2in

\centerline{\large\bf  Harold Rosenberg \, and \, Renato Tribuzy${}^*$}

\footnote{2010 Math. Subject Classification: Primay: 53C42,
Secondary: 53C40.} \footnote{Key words: Rigidity, Bonnet problem,
convex surfaces, homogeneous space.} \footnote{${}^*$ Partially
supported by CNPq and FAPEAM.}

\vskip .5in

\begin{abstract}
   We prove rigidity of oriented isometric immersions of complete surfaces in the homogeneous 3-manifolds  $E(k,\tau)$ (different from the space forms) having the same positive extrinsic curvature.
\end{abstract}

%\vglue .1in

%\noindent{\bf Key words:} Rigidity, Bonnet problem, convex surfaces, homogeneous spaces.

%\vglue .1in

%\noindent{\bf Subject Classification:} Primary 53C42, Secundary 53C40.

\vglue .3in

\noindent{\large\bf Introduction.}

\medskip

An isometric immersion $f\colon M \to N$ is rigid if given any other isometric immersion $g\colon M \to N$, there is an isometry $h\colon N \to N$ such that $hf = g$. An isometric immersion $f\colon M \to N$ is locally rigid if whenever $f(t)\colon M \to N$ is a smooth family of isometric immersions with $f(0)=f$, then there are isometries $h(t)\colon N \to N$ such that $h(t)f(t) = f$.

\smallskip

Strictly convex compact surfaces in $\mathbb{R}^3$ are rigid [C], and there are complete strictly convex surfaces in $\mathbb{R}^3$ that are not rigid [O], [P].  A beautiful open problem is to decide if there is a smooth closed surface $M$ in $\mathbb{R}^3$ that is not locally rigid; i.e., is there a continuous one parameter family of isometric immersions of $M$ into $\mathbb{R}^3$ that are not congruent?

\medskip

In this paper we consider local rigidity of convex surfaces in the 3-dimensional simply connected homogeneous 3-manifolds $E(k,\tau)$, $k-4\tau^2 \ne 0$. After the space forms (isometry group of dimension 6), they are the most symmetric 3-manifolds (isometry group of dimension 4). $E(k,\tau)$ is a Riemannian submersion over the two dimensional space form $M^2(k)$, of curvature $k$: $M^2(k) = S^2(k)$ if $k > 0$, $\mathbb{R}^2$ if $k = 0$, $\mathbb{H}^2(k)$ if $k < 0$. The bundle curvature is $\tau$, and the unit tangent field to the fiber $\xi$ is a Killing field. There is a 2-dimensional group of horizontal translations, translation along the $\xi$-orbits are isometries, and rotations about any vertical fiber. An important discrete group of isometries is generated by rotation by $\pi$ about any horizontal geodesic.

When $\tau = 0$, $E(k,0) = S^2(k) \times \mathbb{R}$ if $k > 0$ and $E(k,0) = \mathbb{H}^2(k) \times \mathbb{R}$, $k < 0$. For $\tau \ne 0$, $k > 0$, they are the Berger spheres. For $\tau \ne 0$, $k = 0$, this gives Nil(3), i.e., Heisenberg space. And $\tau \ne 0$, $k < 0$, $E(k,\tau) = \widetilde{\text{PSL}(2,\mathbb R)}$; the universal covering space of the unit tangent bundle of $\mathbb{H}^2(k)$.

Convexity in $E = {E}(k,\tau)$ can be defined in terms of the second fundamental form. The least one needs is the extrinsic curvature $K_e$ (the product of the principal curvatures) should be positive. However when $\tau \ne 0$, one needs the principal curvatures to be at least $|\tau|$  to obtain global theorems. Assuming $K_e > 0$, one has a Hadamard-Stoker theorem in $\mathbb{H}^2 \times \mathbb{R}$: if $f\colon M^2 \to \mathbb{H}^2\times\mathbb{R}$ is an immersion (complete) with $K_e > 0$, then $f$ is an embedding and $M^2 = S^2$ or $\mathbb{R}^2$. Also one can describe the embedding, [E-G-R].  This theorem is also true in $E(k,\tau)$ provided the principal curvatures are greater than $|\tau|$, [E-R].

In this paper we prove local rigidity of complete surfaces in $E(k,\tau)$ with the same positive extrinsic curvature, and satisfying a three point condition.
 In $E(k,\tau)$, $K$ and $K_e$ are related by the Gauss equation so knowing both these functions tells us the angle (up to sign) the tangent plane of the surface makes with the vertical fiber $\xi$.

In a beautiful paper [G-M-M], the authors studied rigidity of surfaces in $E(k,\tau)$ having the same principal curvatures (the Bonnet problem).

They showed that such surfaces, that are also real analytic, are rigid with some exceptions. In $E(k,0)$ the exceptions are minimal surfaces (they have a 1-parameter family of isometric deformations; like the associated family from the catenoid to the helicoid in $\mathbb{R}^3$) and screw motion helicoidal surfaces. When $\tau \ne 0$, the exceptions are the helicoidal surfaces and Benoit Daniels' CMC twin surfaces [B].

\bigskip

\noindent{\bf The Main Result.}

\medskip

Let $M$ be a complete Riemannian oriented surface. Given an immersion $f\colon M \to E(k,\tau)$ and an oriented frame $(e_1,e_2)$ of $M$, we define the unit normal $N_f$ so that $(f_*(e_1), f_*(e_2), N_f)$ is positive in $E(k,\tau)$ (assumed oriented). When $K_e(f) > 0$, the principal curvatures of $f(M)$ have the same sign on $M$ (always positive or always negative). We always choose the orientation so that they are positive. Given two such immersions with the same positive $K_e$ we say $fg^{-1}\colon g(M) \to f(M)$ is positive when  $fg^{-1}$ is orientation preserving, i.e., both surfaces have positive principal curvatures. We say that $f$ is strictly convex in $E(k,\tau)$ if $K_e(f) > \tau^2$.

\bigskip

\noindent{\bf Theorem A.} Let $f(t)\colon M \to  E = {E}(k,\tau) $ be a smooth family of isometric immersions with $f(0)=f$. Suppose $f$ is strictly convex, $K_e(f_t(x)) = K_e(f(x))$ for $x \in M$ and all $t$, and $Hf_t(x) = H(f(x))$ at three distinct points $x$ of $M$. Then there are isometries $h(t) \colon E \to E$ such that $h(t)f(t)=f$.

\medskip

We begin with a remark on vertical points of the immersion.

\bigskip

\noindent{\bf Lemma 1.} Let $f\colon M \to E$ be an immersion with $K_e(x) > 0$\, $\forall\,x \in M$. Let $g\colon M \to \mathbb{R}$ be the ``angle'' function: $g(x) = \langle N(x),\xi\rangle$, $N$ the unit normal $N_f$\,. If $p \in \Sigma = f(M)$ and $T_p(\Sigma)$ is vertical (i.e. $\xi(p) \in T_p(\Sigma))$, then $f$ is a submersion in a neighborhood of $p$.

\bigskip

\noindent{\bf Corollary 1.} At a vertical point $p \in \Sigma$, there is a disk neighborhood $D$ of $p$ in $\Sigma$ such that $g^{-1}(0)$ is a smooth curve $\beta$ through $p$; $\beta$ separates $D$ into 2 components and $g$ has opposite signs on the two components.

\bigskip

\noindent{\bf Proof of Lemma.}

Let $\pi\colon E(k,\tau) \to M^2(k)$ be the Riemannian submersion, and let $\gamma$ be a geodesic of $M^2(k)$, $\gamma(0) = \pi(p)$, and $d\pi(N(p)) = \gamma'(0)$. Let $P$ be the vertical ``plane'': $P = \pi^{-1}(\gamma)$; $P$ is isometric to $\mathbb{R}^2$ and totally geodesic when $\tau = 0$ (the extrinsic curvature of $P$ is $-\tau^2$). At $p$, $N(p)$ is tangent to $P$, so $\Sigma \cap P$ is a smooth curve $C(s)$, for $s$ near 0; $C(0) = p$, $C'(0) = \xi(p)$.

Since $C$ is a normal section of $\Sigma$ at $p$, the curvature of $C$ at $p$ in $P$ is between the two principal curvatures of $\Sigma$ at $p$, so $k_d^P(p) > 0$. Let $T(s) = C'(s)$, $s$ arc length along $C$.

Denote by $N(s)$ the normal to $\Sigma$ at $C(s)$ and $N_C^P(s)$ the unit normal to $C(s)$ in $P$. Write
$$
N(s) = a(s)\,N^P(s) + E(s),
$$
where $E(s)$ is normal to $P$ along $C(s)$. We have $N(0) = N(p) = N_C^P(0)$ so $a(0) = 1$, $E(0) = 0$.

We want $dg_p(\xi) \ne 0$; we calculate
$$
dg_p(\xi) = \frac{d}{ds}\bigg|_{s=0}\langle \xi,N(s)\rangle = \langle \widetilde{\nabla}_T\xi,N\rangle(0) + \langle\xi,\widetilde{\nabla}_TN\rangle(0)\\
 = \langle \xi,\widetilde{\nabla}_T N\rangle(0).
$$
Since
$$
 \widetilde{\nabla}_T\xi(0) \\= \tau(T\wedge\xi)(0)\\
= \tau(\xi \wedge\xi)(0) = 0; (T(0) = \xi).
$$

Now
$$
\widetilde{\nabla}_TN(s) = a'(s)N_C^P(s) + a(s)\widetilde{\nabla}_T N_C^P(s) + \widetilde{\nabla}_T E(x).
$$
Hence
$$
\langle\widetilde{\nabla}_TN, \xi\rangle(0) = \langle\widetilde{\nabla}_TN_C^P,\xi\rangle(0) + \langle\widetilde{\nabla}_T E,\xi\rangle(0).
$$
Since $E$ is normal to $P$ and $\xi$ is tangent to $P$, $\langle E(s),\xi\rangle = 0$ along $C$. So
$$
\langle\widetilde{\nabla}_TE,\xi\rangle = -\langle E,\widetilde{\nabla}_T\xi\rangle.
$$
At $s=0$, $E(0)=0$, so $\langle\widetilde{\nabla}_TE,\xi\rangle(0) = 0$.

Finally we have
$$
dg_p(\xi) = \langle\xi,\widetilde{\nabla}_T N_C^P\rangle(0).
$$
By the Gauss equation for $P$:
$$
\widetilde{\nabla}_TN_C^P(0) = \nabla_T^P\,N_C^P(0) + \text{II}^P(\xi, N(0))(N(0)\wedge\xi).
$$
Hence
$$
\langle\widetilde{\nabla}_T N_C^P(0),\xi\rangle = \langle \nabla_T^P N_C^P,\xi\rangle(0) = k_C^P(0) \ne 0;
$$
i.e. \, $dg_p(\xi) = k_C^P(0) \ne 0$.

\bigskip

\noindent\textbf{Proof of Theorem A.} Let $f\colon M \to E(k,\tau)$ be a strictly convex isommetric immersion. We will define a special set of moving frames away from the horizontal points $(\xi\,\perp\,TM)$ given by $\varepsilon_1 = \frac{P(\varepsilon)}{|P(\varepsilon)|}$\,, where $P$ denotes the projection into the tangent plane, $J$ is the positive rotation and $\varepsilon_1 = J\varepsilon_2$\,. This way, we can write
$$
\xi = cos(\theta)\varepsilon_1 + sin(\theta)N
$$
where $\theta$ is the function meauring the angle between the vectors $\xi$ and $TM$.  The function $\theta$ and its derivative are defined at least locally.

In order to calculate the second fundamental form, we differentiate $\langle\varepsilon_1,\xi\rangle = cos(\theta)$, getting
$$
X\langle\varepsilon_1,\xi\rangle = \langle\nabla_X\varepsilon_1,\xi\rangle + \langle\alpha(\varepsilon_1,X)N,\xi\rangle + \tau\langle\varepsilon_1,X \wedge \xi) = -sin(\theta)d\theta X.
$$
Because of $\langle\varepsilon_1,\xi\rangle=0$, we also get
$$
\alpha(\varepsilon_1,X)sin(\theta) + \tau(X,\varepsilon_2\rangle sin(\theta) = -sin(\theta)d\theta X
$$
onsequently we have for those points where $sin(\theta) \ne 0$:
\begin{equation} \alpha(\varepsilon_1,X) = -d\theta X - \tau\langle X,\varepsilon_2\rangle. \tag{1}
\end{equation}
Now, for the points where $sin(\theta)=0$, we know they lie in a differential curve as shown in Lemma 1.

By continuity, this equation holds for all points whenever the special moving frame is defined, so this holds for all non-horizontal points.

From (1), we obtain
$$
\alpha(\varepsilon_1,\varepsilon_1) = -d\theta\cdot \varepsilon_1
$$
$$
\alpha(\varepsilon_1,\varepsilon_2) = -d\theta\cdot \varepsilon_2 - \tau
$$
Since the immersion is convex, $\alpha(\varepsilon_1,\varepsilon_1) \ne 0$ and thus $d\theta \ne 0$ at each non-horizontal point.

\bigskip

\noindent\textbf{Lemma 2.} Every horizontal point of a convex immersion, is isolated.

\medskip

\noindent\textbf{Proof:} The horizontal points are the vanishing points of the field $\xi \in N$. Let $p \in M$ be a point so that $(\xi \times N)_p = 0$ and let $\{v_1,v_2\}$ be an orthonormal positively oriented basis which diagonalizes the Weingarten operator at the point $p$ and is associated to the eigenvalues $\lambda_1$ and $\lambda_2$\,.

Let $\overline\nabla$ be the connection of $E(k,\tau)$, then:
\begin{align*}
\overline{\nabla}_{v_1}(\xi \times N)_p &= (\overline{\nabla}_{v_1}\xi \times N)_p + (\xi \times \overline{\nabla}_{v_1}N)_p\\
&= \tau(v_1\times \xi)_p \times N_p - \lambda_1(\xi \times v_1)_p = \tau v_1 - \lambda_1 v_2
\end{align*}
and we also have:
$$
\overline{\nabla}_{v_2}(\xi \times N)_p = \tau(v_2 \times \xi)_p \times N_p - \lambda_2(\xi \times v_2)_p = -\tau v_2 + \lambda_2 v_1\,.
$$

Let $\{V_1,V_2\}$ be a parallel extension of $\{X_1,X_2\}$ along the geodesic which starts at $p$. Let us consider the smooth map $F$ defined in a neighborhood of $p$ in $M$ taking values in $\mathbb{R}^2$ given by
$$
F(g) = (\langle V_1,\xi \times N\rangle_q\,,, \langle V_2,\xi \times N\rangle_q).
$$
Considering that $F'(p)\cdot v_1 = (\tau,-\lambda_1)$ and also $F'(p)\cdot v_2 = (\lambda_2,-\tau)$, we have that $F'(p)$ is an isomorphism.

\medskip

Therefore $F$ is a diffeomorphism when restricted to a neighborhood of $p$. Thus $p$ is the unique zero of the function $F$ and consequently the unique zero of $\xi \times N$.

\bigskip

\noindent\textbf{Corollary 2.} When $M$ is  compact, the set of horizontal points is finite and equals two.

\medskip

\noindent\textbf{Proof:} At every horizontal point we have $cos(\theta)=0$ and thus either $cos(\theta) \ge 0$ or $cos(\theta) \le 0$. Consequently we can choose a particular range for the function $\theta$ either in $\big[-\frac{\pi}{2}, \frac{\pi}{2}\big]$ or $\big[\frac{\pi}{2}, \frac{3\pi}{2}\big]$.

Let us proceed with the calculation of the second fundamental form. Differentiating the equation $\langle \varepsilon,\xi\rangle = 0$, we obtain
\begin{align*}
0 &= X\langle\varepsilon_2,\xi\rangle = \langle\nabla_X \varepsilon_2,\xi\rangle + \langle\alpha(X,\varepsilon_2)N,\xi\rangle + \tau\langle\varepsilon_2,X \times \xi\rangle\\
&= \langle{\nabla}_X \varepsilon_2,\varepsilon_1\rangle cos(\theta) + \alpha(X,\varepsilon_2) sin(\theta) + \tau\langle \varepsilon_1,X\rangle sin(\theta).
\end{align*}
Consequently $\alpha(X,\varepsilon_2) = cotg(\theta)w_{12}(X) - \tau\langle \varepsilon_1,X\rangle$ and it follows that
\begin{align*}
\alpha(\varepsilon_1,\varepsilon_2) &= cotg(\theta)w_{12}(\varepsilon_1)-\tau\\
\alpha(\varepsilon_2,\varepsilon_2) &= cotg(\theta)w_{12}(\varepsilon_2),
\end{align*}
where $w_{ij}(X) = \langle\nabla_X \varepsilon_i,\varepsilon_j\rangle$.

Now, we will determine ordinary differential equations satisfied by a certain angle function we now define.

Let $v \in \mathfrak{X}(M)$ be a unit vector such that $d\theta\cdot v = 0$ and chosen in such a way that $Jv = \frac{grad(v)}{|grad(v)|}\,\cdot$

Let $\phi$ be the angle between the vectors $\varepsilon_1$ and $v$. That is,
\begin{align*}
v &= cos(\phi)\varepsilon_1 + sin(\phi)\varepsilon_2\\
Jv &= -sin(\phi)\varepsilon_1 + cos(\phi)\varepsilon_2\,.
\end{align*}
As $\alpha(v,\varepsilon_1) = -d\theta v - \tau\langle v,\varepsilon_2\rangle$, we have $\alpha(v,\varepsilon_1) = -\tau sin(\phi)$ and consequently
\begin{equation}
cos(\phi)\alpha(\varepsilon_1,\varepsilon_2) + sin(\phi)\alpha(\varepsilon_1,\varepsilon_2) = -\tau sin(\phi). \tag{2}
\end{equation}
In a similar way, $\alpha(Jv,\varepsilon_1) = -d\theta Jv - \tau\langle Jv,\varepsilon_2\rangle = -|grad \theta| - \tau cos(\phi)$ and therefore
\begin{equation}
-sin(\phi)\alpha(\varepsilon_1,\varepsilon_2) + cos(\phi)\alpha(\varepsilon_1,\varepsilon_2) = -|grad \theta| - \tau cos(\phi). \tag{3}
\end{equation}
From equations (2) and (3), we get
\begin{align*}
\alpha_{11} &= |grad \theta| sin(\phi)\\
\alpha_{12} &= -|grad \theta| cos(\phi) - \tau.
\end{align*}
In order to obtain the first differential equation satisfied by $\phi$, we will calculate $\alpha(\varepsilon_2,v)$ using two different approches.

On the one hand we have:
$$
\alpha(\varepsilon_2,v) = cos(\phi)\alpha(\varepsilon_1,\varepsilon_2) + sin(\phi)\alpha(\varepsilon_2\varepsilon_2).
$$
Denoting $K_\varepsilon$ as the extrinsic curvature of $M$, we have
$$
\alpha(\varepsilon_2,\varepsilon_2) = \frac{K_\varepsilon - \alpha(\varepsilon_1,\varepsilon_2)^2}{\alpha(\varepsilon_1,\varepsilon_1)}\,\cdot
$$
Note that since $M$ is convex, $\alpha(\varepsilon_1,\varepsilon_1) \ne 0$, and therefore
\begin{align*}
\alpha(\varepsilon_2,v) &= -cos(\phi)(|grad \theta| cos(\phi)+\tau)\\
&+ sin(\phi)\,\frac{K_\varepsilon - (|grad \theta| cos(\phi)-\tau)}{|grad \theta| sin(\phi)}\,\cdot
\end{align*}
Note that $sin(\phi) \ne 0$ because of $\alpha(\varepsilon_1\varepsilon_1) \ne 0$. Consequently,
\begin{align*}
\alpha(\varepsilon_2,v) &= cos(\phi)(|grad \theta| cos(\phi)+\tau)\\
&+ \frac{K_\varepsilon - (|grad \theta| cos(\phi)-\tau)^2}{|grad \theta|} \cdot \tag{4}
\end{align*}
On the other hand, we have $w_{12} = \widetilde{w}_{12} - d\phi$, where $\widetilde{w}_{12} X = \langle\nabla_X v,Jv\rangle$ and considering that $\alpha(v,\varepsilon_2) = cotg(\theta)w_{12}(v) - \tau\langle\varepsilon_1,v\rangle
$. We have
\begin{equation}
\alpha(v,\varepsilon_2) = cotg(\theta)(\widetilde{w}_{12}(v) - d\phi,v) - \tau cos(\phi). \tag{5}
\end{equation}
Equating the equations (4) and (5), we get a differential equaation satisfied by $\phi$ along the trajectories of $v$. Namely:
\begin{align*}
d\phi\cdot v &= \widetilde{w}_{12}(v) - tg(\theta)\{\tau cos(\phi) + cos(\phi)\cdot[|grad(\theta)|cos(\phi)+\tau]\\
&+ \frac{1}{|grad(\theta)|}\cdot [K_\varepsilon - (|grad(\theta)| cos(\phi)-\tau)^2]\}.
\end{align*}
Similarly, we have $\alpha(Jv,\varepsilon_2) = -sin(\phi)\alpha(\varepsilon_1,\varepsilon_2) + cos(\phi)\alpha(\varepsilon_2,\varepsilon_2)$ and therefore
\begin{align*}
\alpha(Jv,\varepsilon_2) &= sin(\phi)(|grad \theta| cos(\phi)+\tau)\\
&+ cotg\,\phi\,\frac{K_\varepsilon - (|grad \theta|cos \phi + \tau)^2}{|grad \theta|}\,\cdot
\end{align*}
The following equation also holds:
\begin{equation}
\alpha(Jv,\varepsilon_2) = cotg(\theta)(\widetilde{w}_{12}(Jv) - d\phi\cdot Jv) + \tau sin(\phi). \tag{7}
\end{equation}
Now we can equate also the equations (6) and (7) to obtain a differential equation satisfied by $\phi$ along the trajectories of $Jv$.
\begin{align*}
d\phi\cdot Jv &= \widetilde{w}_{12}(Jv) - tg(\theta)\{sin(phi)\cdot[|grad(\theta)| cos(\phi)+\tau]\\
&= cotg(\phi)\,\frac{K_\varepsilon - (|grad(\theta)| cos(\phi)-\tau)^2}{|grad(\theta)|} + \tau sin(\phi)\}.
\end{align*}
Therefore, if two isometric convex immersions have the same extrinsic curvature, the same function $\theta$ and the same function $\phi$ at a point, then they have the same function $\phi$ in a neighborhood of that point.

\medskip

To complete the proof of the Theorem observe:

\medskip

\noindent\textbf{Fact 1:} From the Gauss equation, the set of possibilities for $\theta$ is discrete and therefore $\theta$ is constant along the deformation.

\medskip

\noindent\textbf{Fact 2:} As the set of critical points of the function $\theta$ consists of two points, there exists a point $p \in M$ where $d\theta \ne 0$ and the mean curvature is preserved along the deformation. In a neighborhood of this point, the special moving frame $\{\varepsilon_1(t),\varepsilon_2(t)\}$ is defined for every value of $t$ of the deformation, as well as the functions $\phi_t$\,.

\medskip

\noindent\textbf{Fact 3:} As $H$ and $d\theta$ are both different from zero at $p$, there exists at most two possible values for $\phi_t(p)$.

Indeed, we have:
$$
\alpha(\varepsilon_1(t),\varepsilon_1(t))\{2H-\alpha(\varepsilon_1(t),\varepsilon_1(t))\} - \alpha(\varepsilon_1(t),\varepsilon_2(t))^2 = K_\varepsilon\,.
$$
Since $\alpha(\varepsilon_1(t),\varepsilon_1(t)) = |grad(\theta)| sin(\phi_t)$ and $\alpha(\varepsilon_1(t),\varepsilon_2(t)) = -|grad(\theta)| cos(\phi_t)-\tau$, we obtain
$$
2H|grad(\theta)| sin(\phi_t) + 2\tau|grad(\theta)| cos(\phi_t)-|grad(\theta)|^2 - \tau^2-K_\varepsilon = 0.
$$
We observe an equation of the type $Asin(\phi_t) + Bsin(\phi_t) + C = 0$ yields a quadratic polynomial equation in the variable $cos(\phi_t)$ and consequently has at most two roots, unless all its coefficients are zero, i.e., $A = B = C = 0$.

Therefore $\phi$ is constant along the deformation at the point $p$. Consequently $\phi$ is constant along the deformation in a neighborhood of $p$. Since $\phi$ is preserved, we have that the second fundamental form is preserved and in particular $H$ is preserved.

Let $p_1$ and $p_2$ be points in $M$ where $d\theta = 0$ and let $U = M-\{p_1,p_2\}$. As $M$ is  connected  $U$ is also  connected.

Let $X = \{q \in U/H_t(q) = H_0(q)$ and $\phi_t(q) = \phi_0(q)\}$, where $H_0$ is the mean curvature of $f$ and $\phi_0$ is the function $\phi$ corresponding to $f$.

\medskip

As we observed above, $X$ is an open set $U$. On the other hand the set $X$ is closed since it is intersection of closed sets of $U$. From the hypothesis of the theorem and the observation above, $p \in X$ and thus $X = U$. Consequently, all immersions $f_t$ of the deformation, have the same function $\theta$, the same second fundamental form and the same horizontal directions in $X$. Thus, they are equal to $f$ up to some isometry of $E(k,\tau)$. This conclusion can be extended clearly to $M$.

\newpage

\noindent\textbf{References}

\begin{itemize}
\item[{[C]}] Cohn-Vossen E., \ \textit{Zwei Satze uber die Starrheit der Eilflchen}, \ Nachr. Ges. Will. Gottingen (1927), 125-134.
\item[{[O]}] Olovisnishni-Koff, S., \ \textit{On the bending of infinite convex surfaces}, \ N.S. 18 (60) (1946), 429-440.
\item[{[P]}] Pogorelov, \ \textit{On the rigidity of general infinite convex surfaces with integral curvature $2\pi$}, \ Dokl. Akad. Nauk SSSR (N.S.) 106  (1956), 19-20.
\item[{[G-M-M]}] Gálvez, J., Martínez, A. and Mira, P., \ \textit{The Bonnet problem for Surfaces in homogeneous $3$-manifolds}, \ Communications in Analysis and Geometry, vol. 16, 5, (2008), 907-935.
\item[{[B]}] Benoit Daniel, \ \textit{Isometric Immersions into $3$-dimensional homogeneous manifolds}, \ Comment. Math. Helv. 82 (2007), 87-131.
\item[{[E-G-R]}] Espinar, J., Gálvez, J. and Rosenberg, H., \ \textit{Complete surfaces of positive extrinsic curvature in product spaces}, \ Comment. Math. Helv. 84 (1009), 351-389.
\item[{[E-R]}] Espinar, J., Rosenberg, H., \ \textit{Convex surfaces immersed in homogeneous $3$-manifolds}, preprint.
\end{itemize}

\vglue .3in

\noindent -- Harold Rosenberg, IMPA, Rio de Janeiro, Brazil

\noindent rosen@impa.br

\smallskip

\noindent -- Renato Tribuzy, UFAM, Manaus, AM, Brazil

\noindent tribuzy@pq.cnpq.br

\end{document}